\documentclass{article}

\usepackage[english]{babel}

\usepackage[letterpaper,top=2cm,bottom=2cm,left=3cm,right=3cm,marginparwidth=1.75cm]{geometry}

\usepackage{amsmath}
\usepackage{graphicx}
\usepackage[colorlinks=true, allcolors=blue]{hyperref}
\usepackage{float}
\usepackage{placeins}
\usepackage{multirow}

\usepackage{tikz}
\usetikzlibrary{shapes,arrows,positioning}
\usepackage{graphicx}
\usepackage{subfig}

\title{STM Image Analysis using Autoencoders}
\author{\normalsize Peter Binev, Joshua Moorehead, Ayush Parambath, Luke Parrella, Rori Pumphrey, Miruna Savu 
\thanks{This research was supported by the NSF Grant DMS 2038080 }
}
\date{July 2024}

\begin{document}
\maketitle

\begin{abstract}
This study explores the application of Convolutional Autoencoders (CAEs) for analyzing and reconstructing Scanning Tunneling Microscopy (STM) images of various crystalline lattice structures. We developed two distinct CAE architectures to process simulated STM images of simple cubic, body-centered cubic (BCC), face-centered cubic (FCC), and hexagonal lattices. Our models were trained on $17\times17$ pixel patches extracted from $256\times256$ simulated STM images, incorporating realistic noise characteristics. We evaluated the models' performance using Mean Squared Error (MSE) and Structural Similarity (SSIM) index, and analyzed the learned latent space representations. The results demonstrate the potential of deep learning techniques in STM image analysis, while also highlighting challenges in latent space interpretability and full image reconstruction. This work lays the foundation for future advancements in automated analysis of atomic-scale imaging data, with potential applications in materials science and nanotechnology.
\end{abstract}

\section{Introduction}

Scanning Tunneling Microscopy (STM) has revolutionized our ability to visualize and manipulate matter at the atomic scale, providing unprecedented insights into the structure and properties of materials. Since its invention by Gerd Binnig and Heinrich Rohrer in 1981 \cite{RevModPhys.59.615}, STM has become an indispensable tool in fields ranging from surface physics to nanotechnology. The technique's ability to produce high-resolution images of conductive surfaces has opened new avenues for studying and engineering materials at the nanoscale.

However, the interpretation and analysis of STM images present significant challenges. The complex nature of the quantum tunneling process, combined with various sources of noise and experimental artifacts, can make it difficult to extract meaningful structural information from raw STM data. Traditional analysis methods often rely heavily on human expertise and can be time-consuming and subjective.

In recent years, the rapid advancement of machine learning techniques, particularly deep learning, has shown promise in addressing these challenges. Autoencoders, a class of neural networks designed to learn efficient data representations, have demonstrated success in various image processing tasks, including dimensionality reduction and denoising \cite{10.1145/2689746.2689747}.

This study explores the application of Convolutional Autoencoders (CAEs) to the analysis of STM images, focusing on four common crystalline lattice structures: simple cubic, body-centered cubic (BCC), face-centered cubic (FCC), and hexagonal. Our goal is to develop a model capable of extracting meaningful features from STM images and reconstructing them with high fidelity, potentially revealing structural information that might be obscured by noise or experimental limitations.

We present two distinct CAE architectures, each designed to balance computational efficiency with feature extraction capabilities. These models are trained on a dataset of simulated STM images, carefully constructed to incorporate realistic noise characteristics and variations in imaging conditions. By working with simulated data, we can systematically evaluate our models' performance across different lattice structures and imaging parameters.

The paper is organized as follows: Section 2 provides an overview of the crystalline structures studied and the STM image simulation process. Section 3 details our CAE architectures and the training methodology. Section 4 presents the results of our experiments, including quantitative performance metrics and qualitative assessments of image reconstruction quality. Finally, we discuss the challenges encountered, potential applications of our approach, and directions for future research.

Through this work, we aim to contribute to the development of more sophisticated, automated tools for STM image analysis, potentially accelerating discoveries in materials science and nanotechnology. By bridging the gap between raw STM data and meaningful structural insights, we hope to enhance our ability to study and manipulate matter at the atomic scale.

\subsection{Scanning Tunneling Microscopy}
{\color{black} Scanning Tunneling Microscopy (STM) is a powerful technique that allows for atomic-scale imaging of surfaces. It was invented by Gerd Binnig and Heinrich Rohrer in 1981, for which they later received the Nobel Prize in Physics \cite{RevModPhys.59.615}. STM operates by measuring the tunneling current between a conductive tip and the sample surface, a process that occurs when the tip is brought very close to the sample at a distance of just a few angstroms. The current is highly sensitive to the distance between the tip and the surface, allowing the microscope to achieve atomic resolution.

STM has the unique ability to manipulate individual atoms, which is crucial for constructing precise nanostructures. As described by Eigler and Schweizer \cite{eigler1990positioning}, "Positioning single atoms with a scanning tunneling microscope" was a groundbreaking demonstration of this capability. Furthermore, STM is particularly sensitive to the local density of states at the atomic level, providing detailed information about the electronic properties of the surface \cite{PhysRevB.42.8841}.
%\cite{chen1990probing}. {\color{red} not sure this is a proper citation -- the article it is pointing to is \cite{PhysRevB.42.8841}}

One of the key features of STM is its ability to image conductive surfaces with extremely high resolution, making it an invaluable tool in fields such as surface physics, materials science, and nanotechnology. By allowing scientists to visualize and manipulate individual atoms, STM has opened up new possibilities in the study and application of nanomaterials.}

\subsection{Autoencoders}
{\color{black} First developed in 1986 by Rumelhart \textit{et al}, Autoencoders sought to solve the problem of "backpropogation without a teacher" (unsupervised learning) in neural networks\cite{6300629}. In other words, they employ a type of deep learning method that does not require labeled data.

Autoencoders are trained to develop two distinct functions: the encoder and the decoder. The encoder serves to map a particular dataset to a lower dimensional latent space, whereas the decoder attempts to reconstruct the input data from the latent space. In doing so, the autoencoder determines an "efficient," or compressed, representation of the input data, learning the most important features from which the original data can still be reasonably obtained. 

Common applications of autoencoders include nonlinear dimensionality reduction \cite{10.1126/science.1127647}, anomaly detection, and denoising \cite{10.1145/2689746.2689747}. The nonlinear property of autoencoders allows its performance to far surpass linear dimensionality reduction techniques like principle component analysis (PCA) while avoiding intense computational complexity required by techniques like kernel PCA.  
}

\section{Crystalline Structures - A Few Examples}
Crystalline structures can be broadly classified into various lattice systems, with hexagonal and cubic lattices being two of the most significant.
Below we give a few examples featuring different crystalline structures of carbon (C) and gold (Au).

\subsection{Cubic Lattices - Diamond and Gold}
{\color{black} Diamond has a cubic lattice structure where each carbon atom is bonded to four others in a tetrahedral arrangement. This structure makes diamond extremely hard and gives it excellent visual properties. Diamonds are used in cutting tools, abrasives, and as gemstones in jewelry due to their unmatched hardness and ability to refract light.

\begin{figure}[H]
    \centering
    \includegraphics[width=0.75\linewidth]{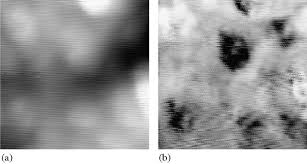}
    \caption{DiamondSTM Image \cite{TROYANOVSKIY2006S27}}
    \label{fig:DiamondSTM}
\end{figure}

Gold forms a face-centered cubic (FCC) lattice. In this structure, each gold atom is surrounded by twelve others, creating a dense and stable arrangement. This structure contributes to gold's high malleability, ductility, and resistance to corrosion. These properties make gold valuable in electronics, jewelry, and as a standard for currency.}

\begin{figure}[H]
    \centering
    \includegraphics[width=0.5\linewidth]{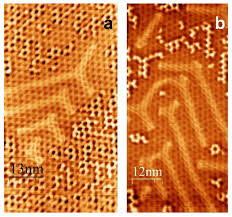}
    \caption{Gold STM Image \cite{ma16103833}}
    \label{fig:GoldSTM}
\end{figure}

\subsection{Hexagonal Lattices - Graphite and Graphene}
{\color{black}Graphite has a hexagonal lattice structure where carbon atoms form layers of hexagons. Each carbon atom is bonded to three others, creating flat sheets that stack on top of each other. These layers can easily slide over one another, which makes graphite a good lubricant and an excellent conductor of electricity. Graphite's structure is useful in applications like pencil leads and as a lubricant in machinery.

Graphene is a single layer of carbon atoms arranged in a hexagonal lattice, similar to one layer of graphite but only one atom thick. Despite its thinness, graphene is incredibly strong—about 200 times stronger than steel—and highly flexible. It also conducts electricity and heat better than most materials. These unique properties make graphene promising for various applications, including electronics, sensors, and even in medical devices.}

\section{Developing STM Data Simulation Software}

\subsection{Simulation Software Overview}
The simulation software creates images with noise to train our autoencoders. To ensure the autoencoders are trained properly the program adds three different types of noise to the STM images: Gaussian, Poisson, and Striation. Additionally, the simulation software also creates images of four different lattice structures: Simple Cubic, Body-Centered Cubic, Face-Centered Cubic, and Hexagonal. 
The setup, the considerations, and the software for the cubic lattices are based on the work of Matthew Belanger \cite{belanger2023}. 
Here we present the specifics for the hexagonal lattices. 

\subsection{Adding Noise to Images}
{\color{black} The simulation program artificially adds Thermal and Shot noise and Striations to the images that it produces. Furthermore, it also adds variations to the position and brightness of the captured atoms. The program refers to Thermal Noise as Gaussian Noise and Shot Noise as Poisson Noise as they follow those distributions respectfully.

Realistically, Gaussian noise is caused by natural sources such as thermal vibration of atoms. The spectral properties of shot, or Poisson, noise is related to the detection of lights reflected from the cantilever \cite{nano11071746}.
For this simulation, separate parameters has been set for Gaussian, Poisson, and Striation noise strength respectively. Using these parameters, a noise matrix is created which the array of points that represent the atoms in the simulated images are passed through. The resulting images are the simulated STM images with various amounts of noise added on. Training the autoencoder on these images strengthens it's ability to denoise and produce clean STM images.}

\subsection{Hexagonal Lattice Structure}
{\color{black}Adding the capability of creating STM images with hexagonal lattice structures greatly improves the capabilities of the autoencoder. Training the autoencoder on these images will advance its ability to denoise STM images of crystals with hexagonal lattice structures such as arsenic, calcite, and quartz \cite{Brit19}.

For the purposes of this paper, we will examine one of two hexagonal lattices. The first step in creating the hexagonal lattice structure is placing a hexagonal grid over the surface of the x,y plane. Each hexagon is made up of separate triangles in order to more simply gather all the relevant atoms for projection. The variable $a$ is used as the length of the base of these triangles. Starting from the origin, points are placed along the $x$ axis at intervals of \textit{a}. The height of each triangle is $a/2$. As we move in the $y$ direction, points are placed at intervals of $3a/2$. Every other row is offset horizontally by $a/2$, this staggering creating the hexagonal grid that we use as the basis for the projection. The atoms are then projected onto the intersection plane after being rotated. Using a randomized tilt angle $\alpha$ and rotation angle, $\theta$ the tilted plane $z$ is created using 

$$ z = \tan(\alpha)\Big(x\cos(\theta) + y\sin(\theta)\Big) $$ 

Note that these trigonometric values have been altered to limit the maximum angle; in this case, the angles \(\theta\) and \(\alpha\) are multiplied by \(\pi/3\), and the angles themselves are found using a \texttt{rand} function between 0 and 1, producing new angles each time the program runs. To determine the appropriate projection parameters, including its brightness and position, for each atom, the simulation uses the vector projection formula. First, a normal vector $\mathbf{N}$ is created perpendicular to the intersection plane which is just a rewritten form of the $z$ plane. 

$$ \mathbf{N} = \big(\tan(\alpha)\cos(\theta) , \tan(\alpha)\sin(\theta) , -1\big)$$ 

We then calculate the vector from each point on the plane by subtracting the coordinates of each point from the origin 

$$ points = (0, 0, 0) - (x, y, \texttt{floor}(z)) $$ 

Each of these point vectors is projected onto the normal vector \(\mathbf{N} = [A, B, C]\) using the projection formula:

\[
\text{proj} = \left(\frac{points \cdot \mathbf{N}}{\mathbf{N} \cdot \mathbf{N}}\right) \cdot \mathbf{N}
\]

The distance of each point from the plane is then calculated by finding the magnitude of the proj vector. Finally, in order to format them for use in testing images, the 3D coordinates are transformed into a 2D plane using a linear projection, where we now have projX, a projected x plane, and projY, a projected y plane. Using a third randomized angle called \(\phi\), a matrix is created and multiplied by the projX and projY to produce the final 2D coordinates in column vector form, shown here as \textit{u} and \textit{v}.

\[
\begin{bmatrix} u \\ v \end{bmatrix} = \begin{bmatrix} \cos(\phi \cdot \pi) & \sin(\phi \cdot \pi) \\ -\sin(\phi \cdot \pi) & \cos(\phi \cdot \pi) \end{bmatrix} \begin{bmatrix} \text{projX} \\ \text{projY} \end{bmatrix}
\]

Before being placed on the image, a variable called spread needs to be applied to each of these points.
}

\subsection{Determining Spread}
{\color{black} Creating clean and accurate STM images through the simulation was a high priority goal. The purpose of the simulation was to mass produce images to efficiently train autoencoders. Therefore, one of the most crucial factors in creating accurate images for each lattice structure was the spread variable. Previously a simple quadratic formula that considered the assigned lattice structure value (1-3 for each lattice structure) was used to determine the spread. Creating the hexagonal lattice structure simulation emphasized that the formula simply was not accurate enough. Through trials, it was determined that there was a unique spread for each lattice structure and a formula would not describe the relationship effectively. Running trials on the simulation established the spreads listed in Table 1 for each of the lattice structures.}
\begin{table}[h]
    \centering
    \begin{tabular}{cccccc}
        Type & Simple Cubic & BCC & FCC & Hexagonal 1 & Hexagonal 2 \\
        Spread & 10 & 13 & 18 & 14 & 10 \\
    \end{tabular}
    \caption{Table outlining the spread for each lattice type}
    \label{tab:my_label}
\end{table}

Each of the coordinates (\textit{u,v}) is now scaled by the spread found above, plotted onto the coordinate grid, and finally the image is cropped. Some example images of the new types, Hexagonal 1 and 2, are shown in the next section.
\subsection{Example Images}
\begin{figure}[H]
    \centering
    % Smaller Type 4 images in one line
    \begin{minipage}{0.3\linewidth}
        \centering
        \includegraphics[width=\linewidth]{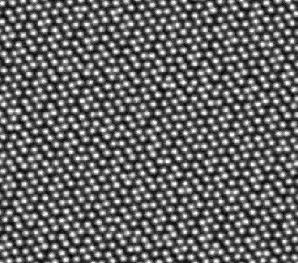}
        \caption{\!Hexagonal \!1, Image \!1}
        \label{fig:label1}
    \end{minipage}
    \hfill
    \begin{minipage}{0.3\linewidth}
        \centering
        \includegraphics[width=\linewidth]{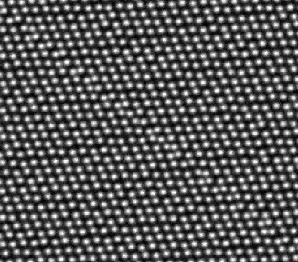}
        \caption{\!Hexagonal \!1, Image \!2}
        \label{fig:label2}
    \end{minipage}
    \hfill
    \begin{minipage}{0.3\linewidth}
        \centering
        \includegraphics[width=\linewidth]{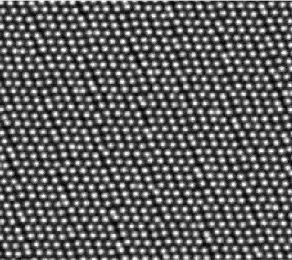}
        \caption{\!Hexagonal \!1, Image \!3}
        \label{fig:label3}
    \end{minipage}
\end{figure}
%    \vspace{1em} 
    
    % Text section between Type 4 and Type 5 images
%    \centering
%    \textbf\\
    Hexagonal 2 is created in a similar way to Hexagonal 1, with slightly different base geometry, and the differences are easily seen when looking at them next to each other. Hexagonal 1 has very close together hexagons, where Hexagonal 2 has clearly defined rows that create spaced out geometry.

%    \vspace{1em} 
    
    % Type 5 images centered in a new row
%    \hspace*{\fill} % Add space before the first image
\begin{figure}[H]    
    \begin{minipage}{0.35\linewidth}
        \centering
        \includegraphics[width=\linewidth]{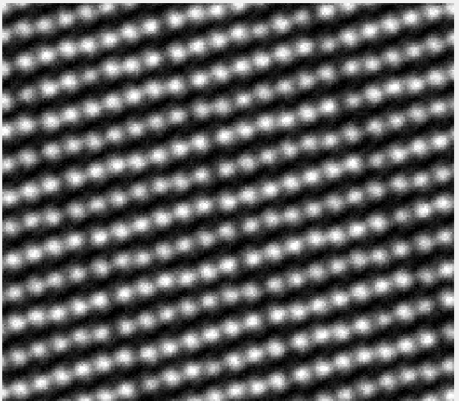}
        \caption{Hexagonal 2, Image 1}
        \label{fig:label4}
    \end{minipage}
    \hspace{2em} % Add space between the two images
    \begin{minipage}{0.35\linewidth}
        \centering
        \includegraphics[width=\linewidth]{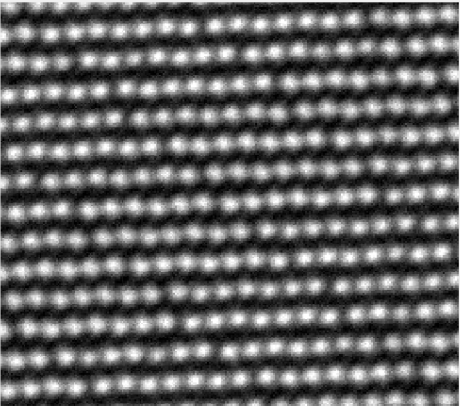}
        \caption{Hexagonal 2, Image 2}
        \label{fig:label5}
    \end{minipage}
    \hspace*{\fill} % Add space after the second image
\end{figure}

\section{Developing Deep Neural Networks Autoencoders for Pattern Extraction from STM Data}

In this section, we present our approach to developing a Convolutional Autoencoder (CAE) for extracting patterns and denoising Scanning Tunneling Microscopy (STM) images. Our goal was to create a model capable of learning compact representations of various crystalline lattice structures while preserving key structural information and ensuring translation invariance. In utilizing an autoencoder to compress simulated images of crystalline lattice structures, the essential features of the structures included in the training dataset (i.e., simple cubic, BCC, FCC, and hexagonal) may be extracted while preserving the structural integrity of the image.

The CAE was designed to process simulated STM images of simple cubic, body-centered cubic (BCC), face-centered cubic (FCC), and hexagonal lattice structures. By training on these diverse structures, we aimed to create a versatile model capable of handling a wide range of crystalline materials. The following subsections detail our methodology, from the architecture design to the evaluation of results, highlighting the challenges we faced and the insights we gained throughout the process.

\subsection{Convolutional Autoencoder Architectures}

In our study, we developed and compared two distinct Convolutional Autoencoder (CAE) architectures, each designed to address the challenge of extracting patterns from STM images. We'll refer to these as CAE-A and CAE-B.

\subsubsection{CAE-A Architecture}

CAE-A, developed by Josh, focuses on efficient feature extraction while maintaining a relatively shallow network structure. This architecture balances model complexity with computational efficiency.

\paragraph{Encoder}
The encoder progressively reduces spatial dimensions while increasing feature depth:
\begin{itemize}
    \item Input layer: 17x17 grayscale image patches
    \item Convolutional layer 1: 16 filters, 3x3 kernel, stride 1, ReLU activation
    \item Max pooling layer 1: 2x2 pool size, stride 2
    \item Convolutional layer 2: 9 filters, 3x3 kernel, stride 1, ReLU activation
    \item Max pooling layer 2: 2x2 pool size, stride 2
    \item Flatten layer
    \item Dense layer: 10 neurons (latent space representation)
\end{itemize}

\paragraph{Decoder}
The decoder reconstructs the original image size from the latent representation:
\begin{itemize}
    \item Dense layer: 144 neurons (to be reshaped to 9x4x4)
    \item Reshape layer: Transforms 144 neurons to 9x4x4 tensor
    \item Transposed convolutional layer 1: 16 filters, 3x3 kernel, stride 2, ReLU activation
    \item Transposed convolutional layer 2: 1 filter, 5x5 kernel, stride 2
\end{itemize}

This architecture captures hierarchical features of the STM images. The latent space dimension of 10 compresses the input while retaining sufficient information for reconstruction. The decoder's first dense layer expands the 10-dimensional latent space to 144 neurons, which are then reshaped to a 9x4x4 tensor. This allows the subsequent transposed convolutions to gradually upsample the feature maps back to the original 17x17 image size.

\begin{figure}[H]
    \centering
    \begin{tikzpicture}[
        node distance=0.4cm,
        auto,
        block/.style={rectangle, draw, text centered, rounded corners, minimum height=1.8em, font=\scriptsize},
        latent/.style={rectangle, draw, fill=lightgray, text centered, rounded corners, minimum height=1.8em, font=\scriptsize},
        line/.style={draw, -latex', shorten >=1pt, shorten <=1pt}
    ]
        
        \node [block, text width=3cm] (input) {Input (17x17x1)};
        \node [block, below=of input, text width=3cm] (conv1) {Conv2D (17x17x16)};
        \node [block, below=of conv1, text width=2.5cm] (pool1) {MaxPool (8x8x16)};
        \node [block, below=of pool1, text width=2.5cm] (conv2) {Conv2D (8x8x9)};
        \node [block, below=of conv2, text width=2cm] (pool2) {MaxPool (4x4x9)};
        \node [block, below=of pool2, text width=2cm] (flatten) {Flatten (144)};
        \node [block, below=of flatten, text width=1.5cm] (dense1) {Dense (10)};
        \node [latent, below=of dense1, text width=1.5cm] (latent) {Latent (10)};
        \node [block, below=of latent, text width=2cm] (dense2) {Dense (144)};
        \node [block, below=of dense2, text width=2cm] (reshape) {Reshape (4x4x9)};
        \node [block, below=of reshape, text width=2.5cm] (tconv1) {TConv2D (8x8x16)};
        \node [block, below=of tconv1, text width=3cm] (tconv2) {TConv2D (17x17x1)};
        \node [block, below=of tconv2, text width=3cm] (output) {Output (17x17x1)};
        
        \path [line] (input) -- (conv1);
        \path [line] (conv1) -- (pool1);
        \path [line] (pool1) -- (conv2);
        \path [line] (conv2) -- (pool2);
        \path [line] (pool2) -- (flatten);
        \path [line] (flatten) -- (dense1);
        \path [line] (dense1) -- (latent);
        \path [line] (latent) -- (dense2);
        \path [line] (dense2) -- (reshape);
        \path [line] (reshape) -- (tconv1);
        \path [line] (tconv1) -- (tconv2);
        \path [line] (tconv2) -- (output);
    \end{tikzpicture}
    \caption{Size-correlated representation of the CAE-A architecture}
    \label{fig:cae_a_architecture}
\end{figure}
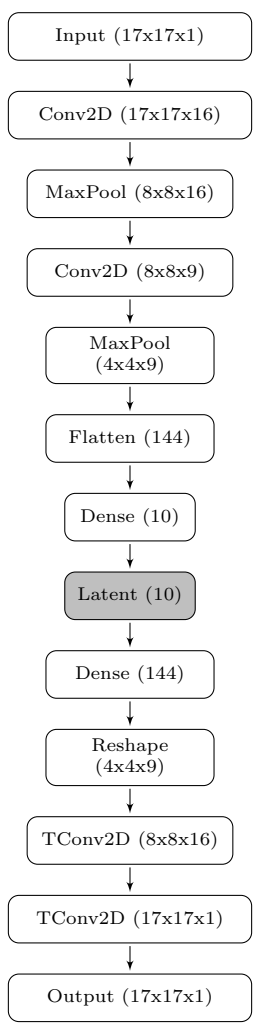

Figure \ref{fig:cae_a_architecture} visualizes the CAE-A architecture, showing data flow through encoder and decoder components.

\subsubsection{CAE-B Architecture}

CAE-B, developed by Rori, takes a different approach to the autoencoder design. While maintaining general structural similarity to CAE-A, CAE-B employs Leaky ReLU activation functions to mitigate risk of "dying neurons" and features more convolutional/transposed convolutional layers in the encoder and decoder, respectively. While training was more computationally expensive for CAE-B, its primary function was to gauge whether increases in model complexity and training time were worth any resulting performance gains. Additionally, CAE-B was developed in MATLAB using the deep learning toolbox.

\paragraph{Encoder}
The encoder of CAE-B is composed of the following layers:
\begin{itemize}
\item Input layer: Image patches of size 16x16 pixels are loaded from an image datastore 
\item Convolutional layer 1: 32 filters, 3x3 kernel, stride of 1
\item Leaky ReLU layer 1: Negative inputs scaled by 0.01
\item Max pooling layer 1: 2x2 pool size, stride of 2
\item Convolutional layer 2: 24 filters, 3x3 kernel, stride of 1
\item Leaky ReLU layer 2: Negative inputs scaled by 0.01
\item Max pooling layer 2: 2x2 pool size, stride of 2
\item Convolutional layer 3: 16 filters, 3x3 kernel, stride of 1
\item Leaky ReLU layer 3: Negative inputs scaled by 0.01
\item Max pooling layer 3: 2x2 pool size, stride of 2
\item Convolutional layer 4: 8 filters, 3x3 kernel, stride of 1
\item Leaky ReLU layer 4: Negative inputs scaled by 0.01
\item Max pooling layer 4: 2x2 pool size, stride of 2
\item Fully connected layer: 10 channels (latent space dimension) 
\end{itemize}

\paragraph{Decoder} 
The decoder of CAE-B is composed of the following layers: 
\begin{itemize}
    \item Project and reshape layer: Projects latent vector into 3D matrix [2 2 24]
    \item Transposed convolutional layer 1: 24 filters, 3x3 kernel, stride of 2
    \item Leaky ReLU layer 1: Negative inputs scaled by 0.01
    \item Transposed convolutional layer 2: 16 filters, 3x3 kernel, stride of 2
    \item Batch normalization layer: Mini-batch data normalized for each independent channel 
    \item Leaky ReLU layer 2: Negative inputs scaled by 0.01
    \item Transposed convolutional layer 3: 8 filters, 3x3 kernel, stride of 2
    \item Leaky ReLU layer 3: Negative inputs scaled by 0.01
    \item Transposed convolutional layer 4: 4 filters, 3x3 kernel, stride of 2
    \item Clipped ReLU layer: Forces output within the range [0,1]
\end{itemize}

\begin{figure}[H]
    \centering
    \begin{tikzpicture}[
        node distance=0.25cm,
        auto,
        block/.style={rectangle, draw, text centered, rounded corners, minimum height=1.3em, font=\tiny},
        latent/.style={rectangle, draw, fill=lightgray, text centered, rounded corners, minimum height=1.3em, font=\tiny},
        line/.style={draw, -latex', shorten >=1pt, shorten <=1pt}
    ]
        
        \node [block, text width=2.5cm] (input) {Input (16x16x1)};
        \node [block, below=of input, text width=2.5cm] (conv1) {Conv2D (16x16x32)};
        \node [block, below=of conv1, text width=2.5cm] (lrelu1) {Leaky ReLU (16x16x32)};
        \node [block, below=of lrelu1, text width=2.2cm] (pool1) {MaxPool (8x8x32)};
        \node [block, below=of pool1, text width=2.2cm] (conv2) {Conv2D (8x8x24)};
        \node [block, below=of conv2, text width=2.2cm] (lrelu2) {Leaky ReLU (8x8x24)};
        \node [block, below=of lrelu2, text width=1.9cm] (pool2) {MaxPool (4x4x24)};
        \node [block, below=of pool2, text width=1.9cm] (conv3) {Conv2D (4x4x16)};
        \node [block, below=of conv3, text width=1.9cm] (lrelu3) {Leaky ReLU (4x4x16)};
        \node [block, below=of lrelu3, text width=1.6cm] (pool3) {MaxPool (2x2x16)};
        \node [block, below=of pool3, text width=1.6cm] (conv4) {Conv2D (2x2x8)};
        \node [block, below=of conv4, text width=1.6cm] (lrelu4) {Leaky ReLU (2x2x8)};
        \node [block, below=of lrelu4, text width=1.3cm] (pool4) {MaxPool (1x1x8)};
        \node [block, below=of pool4, text width=1cm] (fc) {Dense (10)};
        \node [latent, below=of fc, text width=1cm] (latent) {Latent (10)};
        \node [block, below=of latent, text width=1.3cm] (project) {Proj. (2x2x24)};
        \node [block, below=of project, text width=1.6cm] (tconv1) {TConv2D (4x4x24)};
        \node [block, below=of tconv1, text width=1.6cm] (lrelu5) {Leaky ReLU (4x4x24)};
        \node [block, below=of lrelu5, text width=1.9cm] (tconv2) {TConv2D (8x8x16)};
        \node [block, below=of tconv2, text width=1.9cm] (bn) {Batch Norm (8x8x16)};
        \node [block, below=of bn, text width=1.9cm] (lrelu6) {Leaky ReLU (8x8x16)};
        \node [block, below=of lrelu6, text width=2.2cm] (tconv3) {TConv2D (16x16x8)};
        \node [block, below=of tconv3, text width=2.2cm] (lrelu7) {Leaky ReLU (16x16x8)};
        \node [block, below=of lrelu7, text width=2.5cm] (tconv4) {TConv2D (16x16x4)};
        \node [block, below=of tconv4, text width=2.5cm] (crelu) {Clipped ReLU (16x16x4)};
        \node [block, below=of crelu, text width=2.5cm] (output) {Output (16x16x1)};
        
        \foreach \i/\j in {input/conv1,conv1/lrelu1,lrelu1/pool1,pool1/conv2,conv2/lrelu2,lrelu2/pool2,pool2/conv3,conv3/lrelu3,lrelu3/pool3,pool3/conv4,conv4/lrelu4,lrelu4/pool4,pool4/fc,fc/latent,latent/project,project/tconv1,tconv1/lrelu5,lrelu5/tconv2,tconv2/bn,bn/lrelu6,lrelu6/tconv3,tconv3/lrelu7,lrelu7/tconv4,tconv4/crelu,crelu/output}
            \path [line] (\i) -- (\j);
    \end{tikzpicture}
    \caption{Size-correlated representation of the CAE-B architecture}
    \label{fig:cae_b_architecture}
\end{figure}
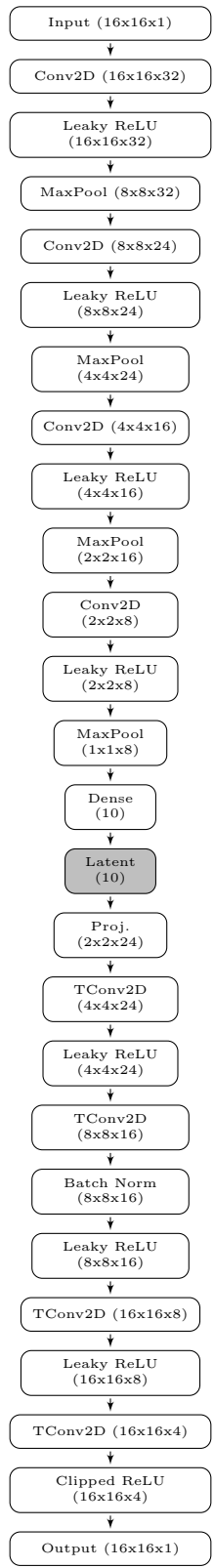

Figure \ref{fig:cae_b_architecture} visualizes the CAE-B architecture, showing data flow through encoder and decoder components.

The development of these two distinct architectures allowed us to explore different approaches to feature extraction and reconstruction in STM images, providing valuable insights into the strengths and limitations of various autoencoder designs for this specific application.

\subsection{Data Preprocessing and Augmentation}

Our data preprocessing and augmentation strategy was designed to prepare the simulated STM images for efficient training of the Convolutional Autoencoder while preserving the realistic noise characteristics inherent in the data.

\subsubsection{Image Preprocessing}

Our preprocessing pipeline consisted of the following steps:

\begin{enumerate}
    \item \textbf{Normalization:} We normalized the pixel values of each image to the range [-1, 1] using the following transformation:
    \begin{equation}
        x_{\mathrm{normalized}} = \frac{2x - (x_{\max} + x_{\min})}{x_{\max} - x_{\min}}
    \end{equation}
    where $x$ is the original pixel value, and $x_{\max}$ and $x_{\min}$ are the maximum and minimum pixel values in the image, respectively. This normalization ensures consistent input scales across all samples and can help with faster convergence during training.
    
    \item \textbf{Patch Extraction:} From each 256x256 STM image, we extracted multiple 17x17 or 16x16 pixel patches. This approach allowed us to focus on local structural features and significantly increased our training dataset size. The patches were extracted with a stride of 4 pixels, resulting in a total of 3600 patches per image.
\end{enumerate}

It's important to note that we did not apply any additional noise reduction techniques, as our simulated STM images already incorporated realistic noise characteristics, including Gaussian and Poisson noise, as well as striations. This approach ensures that our model learns to handle the types of noise typically present in real STM data.

\subsubsection{Data Augmentation}

To enhance the robustness and generalization capability of our model, we implemented the following augmentation techniques:

\begin{enumerate}
    \item \textbf{Random Rotations:} We applied random rotations of 0°, 90°, 180°, or 270° to each patch, enhancing the model's rotational invariance.
    
    \item \textbf{Random Flips:} Horizontal and vertical flips were randomly applied, further increasing the diversity of our training data.
\end{enumerate}

These augmentation techniques were applied on-the-fly during training, effectively expanding our dataset and improving the model's ability to generalize across various crystal orientations.

\begin{figure}[H]
    \centering
    \subfloat[Original 256x256 STM image]{
        \includegraphics[width=0.4\textwidth]{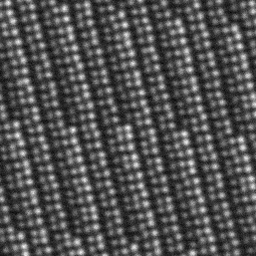}
        \label{fig:full_image}
    }
    \hfill
    \subfloat[Extracted 17x17 patches]{
        \includegraphics[width=0.4\textwidth]{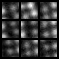}
        \label{fig:patches}
    }
    \caption{Comparison of original STM image and extracted patches}
    \label{fig:image_and_patches}
\end{figure}

As shown in Figure \ref{fig:image_and_patches}, we can see the original STM image (Figure \ref{fig:full_image}) alongside examples of the extracted patches (Figure \ref{fig:patches}).

Our preprocessing and augmentation strategies were designed to preserve the critical structural information and realistic noise characteristics of the simulated STM images while enhancing the model's ability to handle variations in crystal orientations. This approach proved crucial in improving the autoencoder's performance and generalization capabilities across different lattice structures.

\subsection{Training Process and Hyperparameter Tuning}
The training of our Convolutional Autoencoder (CAE) focused on simple cubic lattice structures, involving careful selection of hyperparameters and optimization techniques to achieve the best performance in reconstructing STM images while extracting meaningful latent representations.

\subsubsection{Training Process}
We trained our model using the following setup:
\begin{itemize}
    \item \textbf{Dataset:} Our training set consisted of simulated STM images of simple cubic lattice structures, with 3000 17x17 patches extracted from each 256x256 image.
    \item \textbf{Hardware:} Training was performed on 4 NVIDIA Tesla T4 GPUs in parallel.
    \item \textbf{Framework:} We used PyTorch version 2.3.1 with CUDA 11.8 for GPU acceleration.
    \item \textbf{Optimizer:} We used the Adam optimizer.
    \item \textbf{Loss Function:} Mean Squared Error (MSE) was used as our loss function, measuring the pixel-wise difference between input and reconstructed images.
\end{itemize}

\subsubsection{Hyperparameter Tuning}
We conducted multiple training sessions with varying hyperparameters to optimize our model's performance. The key configurations explored were:
\begin{enumerate}
   \item \textbf{Baseline Configuration:}
   Learning rate = 0.001, Batch size = 1024, Patches per image = 3000, Epochs = 100
   
   \item \textbf{Lower Learning Rate:}
   Learning rate = 0.0001, Batch size = 1024, Patches per image = 3000, Epochs = 100
   
   \item \textbf{Small Batch:}
   Learning rate = 0.001, Batch size = 256, Patches per image = 3000, Epochs = 100
   
   \item \textbf{Large Batch:}
   Learning rate = 0.002, Batch size = 2048, Patches per image = 3000, Epochs = 100
   
   \item \textbf{More Patches:}
   Learning rate = 0.001, Batch size = 1024, Patches per image = 4900, Epochs = 100
   
   \item \textbf{Extended Training:}
   Learning rate = 0.001, Batch size = 1024, Patches per image = 3000, Epochs = 200
   
   \item \textbf{Learning Rate Decay:}
   Initial learning rate = 0.001, Batch size = 1024, Patches per image = 3000, Epochs = 100
\end{enumerate}

To prevent overfitting, we monitored validation loss and implemented comprehensive logging of training metrics. The performance of each configuration was evaluated across three different lattice structures: Simple Cubic, Body-Centered Cubic (BCC), and Face-Centered Cubic (FCC).

\begin{figure}[H]
    \centering
    \includegraphics[width=\textwidth]{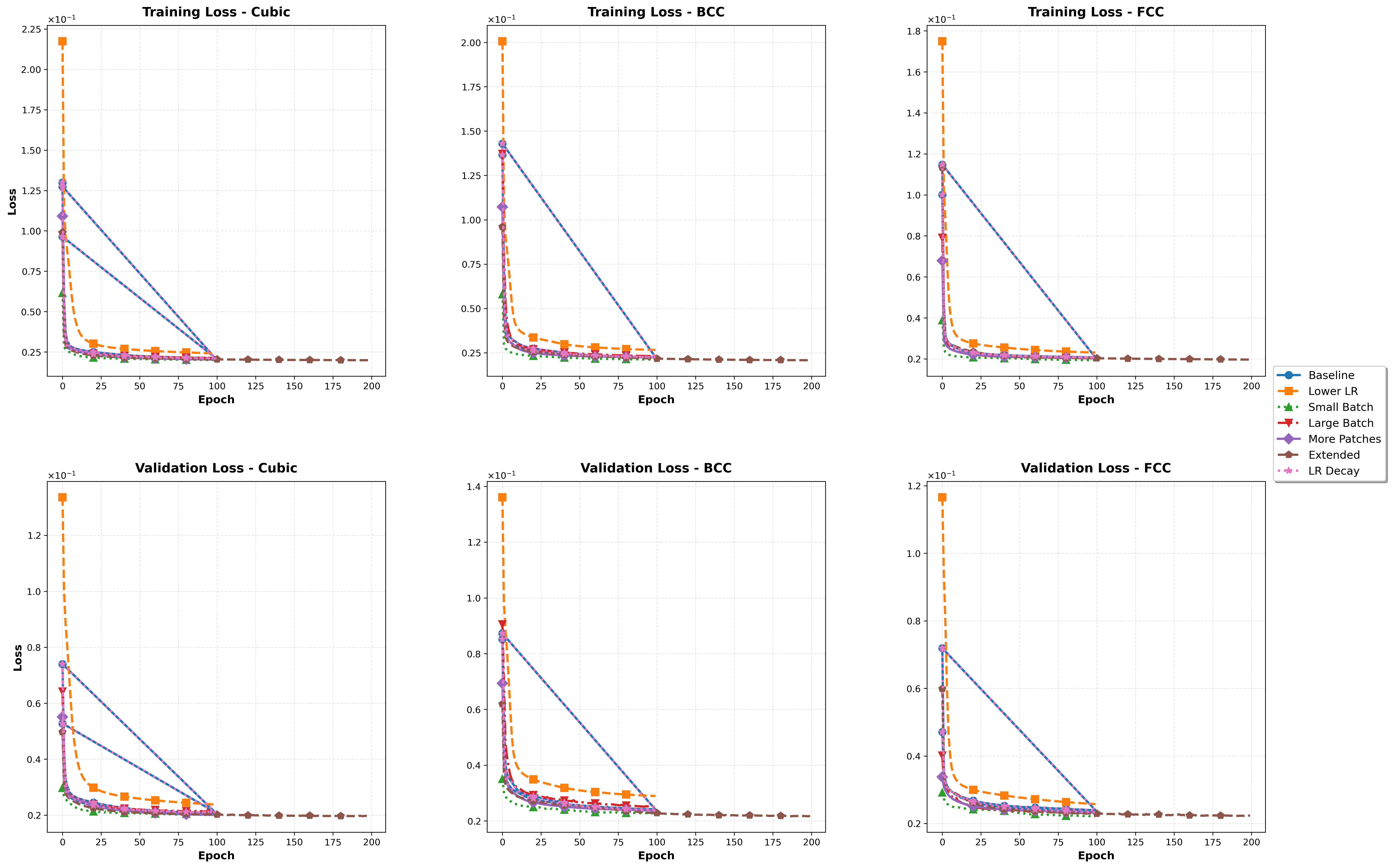}
    \caption{Training (top row) and validation (bottom row) loss curves across different lattice structures for various hyperparameter configurations. Each column represents a different lattice type (Cubic, BCC, and FCC from left to right). The plots reveal distinct learning patterns across configurations and lattice complexities, with notable differences in convergence rates and final performance metrics.}
    \label{fig:loss_curves}
\end{figure}

Analysis of the training and validation curves (Figure \ref{fig:loss_curves}) reveals several key findings across configurations and lattice types:

\begin{itemize}
   \item \textbf{Learning Rate Impact:} The lower learning rate configuration (0.0001, orange dashed line) exhibits consistently higher initial loss values and significantly slower convergence across all lattice types. While it eventually converges, our baseline learning rate of 0.001 (blue solid line) achieves better performance with faster convergence, particularly evident in the rapid initial descent phase.
   
   \item \textbf{Batch Size Effects:} The small batch configuration (256, green dotted line) demonstrates more volatile training curves but achieves comparable or slightly better final performance than larger batches (2048, red dash-dotted line). This is particularly noticeable in the Simple Cubic structure, suggesting that more frequent parameter updates can be advantageous despite increased training noise.
   
   \item \textbf{Extended Training Benefits:} The extended training configuration (brown dashed line) reveals continued improvement beyond the standard 100 epochs, most notably in the FCC structure. The extended training period allows the model to achieve lower final loss values, particularly beneficial for more complex lattice structures.
   
   \item \textbf{Structural Complexity Correlation:} A clear relationship emerges between lattice complexity and convergence behavior. The Simple Cubic structure exhibits the lowest loss values and fastest convergence (visible in both training and validation curves), while FCC shows higher loss values and requires more epochs to stabilize, reflecting its greater structural complexity. BCC demonstrates intermediate behavior, aligning with its moderate structural complexity.
   
   \item \textbf{Generalization Performance:} The close tracking between training (top row) and validation (bottom row) curves across all configurations indicates robust generalization without significant overfitting. The learning rate decay configuration (pink dotted line) particularly demonstrates stable generalization, with minimal divergence between training and validation performance.
   
   \item \textbf{Configuration Stability:} The baseline configuration consistently achieves a good balance between convergence speed and final performance across all lattice types, though its advantage is most pronounced in simpler structures. The more patches configuration (purple diamond line) shows similar stability but with marginally better final performance in some cases.
\end{itemize}

These insights informed our final training strategy: simpler structures like Simple Cubic benefit most from the baseline configuration with standard training duration, while complex structures like FCC show improved results with extended training periods. The choice between small and large batch sizes can be tailored to computational resources, as both achieve comparable final performance despite different convergence characteristics.

Throughout the training process, we used TensorBoard to visualize and monitor these metrics in real-time, enabling rapid identification of optimal configurations for each lattice structure. This comprehensive monitoring approach proved essential for fine-tuning the model's performance across different structural complexities.

\subsection{Latent Space Analysis}

Our Convolutional Autoencoder (CAE) encodes each 17x17 pixel STM image patch into a 10-dimensional latent vector. To visualize this high-dimensional space, we employed Principal Component Analysis (PCA) to project the latent vectors onto a 3D space.

\begin{figure}[H]
    \centering
    \includegraphics[width=0.9\linewidth]{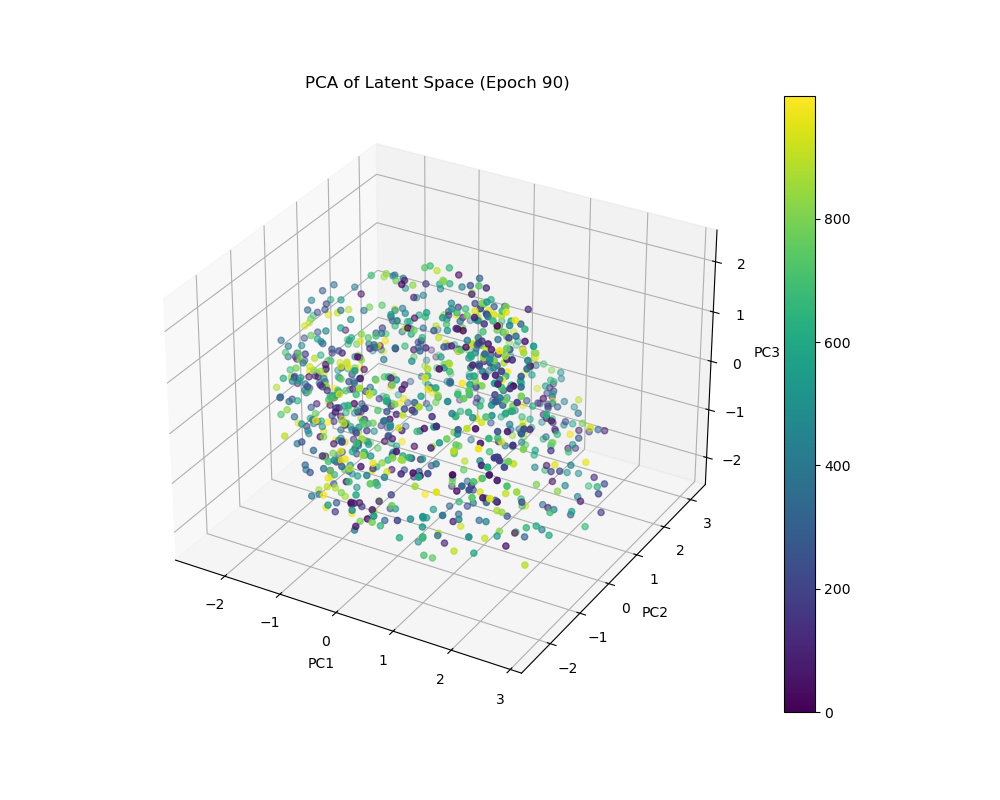}
    \caption{PCA visualization of the latent space at epoch 100}
    \label{fig:pca_latent_space}
\end{figure}

Figure \ref{fig:pca_latent_space} shows the PCA visualization of our latent space after 100 epochs of training. Each point represents a single STM image patch, with colors indicating different values. The visualization reveals a relatively homogeneous distribution of latent vectors, without clear clustering based on lattice structures. This observation suggests potential challenges in our current approach, which are discussed in detail in Section 4.7 (Challenges and Future Work).

\subsection{Results and Performance Evaluation}

We evaluated our Convolutional Autoencoder (CAE) model on simulated STM image patches of various lattice structures. The performance was assessed using both quantitative metrics and qualitative visual inspection at the patch level.

\subsubsection{Quantitative Evaluation}
Our evaluation revealed several key findings across lattice types and configurations:

\begin{itemize}
    \item \textbf{Simple Cubic:} The More Patches configuration achieved the best MSE (0.0152) while Large Batch configuration achieved the best SSIM (0.9270). The Lower LR configuration performed notably worse, suggesting that the default learning rate was more appropriate for this lattice structure.
    
    \item \textbf{BCC:} Extended Training showed the best MSE (0.0185) while Small Batch and More Patches configurations tied for best SSIM (~0.909). This suggests that BCC structures might benefit from longer training periods or increased data sampling.
    
    \item \textbf{FCC:} Extended Training configuration significantly outperformed others, achieving both the best MSE (0.0152) and SSIM (0.9491) across all lattice types and configurations. This indicates that FCC structures might require more training epochs to fully capture their complex patterns.
\end{itemize}

Notably, the Lower LR configuration consistently underperformed across all lattice types, suggesting that our baseline learning rate of 0.001 was more appropriate for this task. The Extended Training configuration showed particular promise for more complex structures like FCC, while simpler structures like Simple Cubic benefited more from increased data sampling (More Patches) or larger batch sizes.

Table \ref{tab:performance_metrics} summarizes the average performance metrics for patches of each lattice type:

\begin{table}[h]
\centering
\begin{tabular}{|c|c|c|c|}
\hline
Lattice Type & Configuration & Average MSE & Average SSIM \\
\hline
\multirow{6}{*}{Cubic} 
& Baseline & 0.0156 & 0.8917 \\
& Lower LR & 0.0326 & 0.8699 \\
& Small Batch & 0.0176 & 0.9130 \\
& Large Batch & 0.0166 & 0.9270 \\
& More Patches & 0.0152 & 0.9027 \\
& Extended Training & 0.0296 & 0.9130 \\
\hline
\multirow{6}{*}{BCC}
& Baseline & 0.0188 & 0.8542 \\
& Lower LR & 0.0292 & 0.7959 \\
& Small Batch & 0.0205 & 0.9093 \\
& Large Batch & 0.0199 & 0.8902 \\
& More Patches & 0.0196 & 0.9092 \\
& Extended Training & 0.0185 & 0.8858 \\
\hline
\multirow{6}{*}{FCC}
& Baseline & 0.0224 & 0.8987 \\
& Lower LR & 0.0261 & 0.8374 \\
& Small Batch & 0.0266 & 0.7656 \\
& Large Batch & 0.0184 & 0.8944 \\
& More Patches & 0.0167 & 0.8768 \\
& Extended Training & \textbf{0.0152} & \textbf{0.9491} \\
\hline
\end{tabular}
\caption{Performance metrics across different lattice types and configurations. Best overall MSE and SSIM values are highlighted in bold.}
\label{tab:performance_metrics}
\end{table}

\subsubsection{Qualitative Evaluation}
Visual inspection of the reconstructed patches provides insight into the model's ability to capture and reproduce the distinctive features of each lattice type at a local level.

\begin{figure}[H]
   \centering
   \begin{minipage}{0.32\linewidth}
       \centering
       \includegraphics[width=\linewidth]{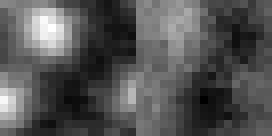}
       \includegraphics[width=\linewidth]{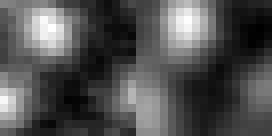}
       {Simple Cubic}
   \end{minipage}
   \begin{minipage}{0.32\linewidth}
       \centering
       \includegraphics[width=\linewidth]{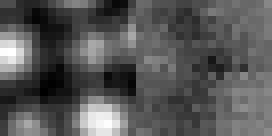}
       \includegraphics[width=\linewidth]{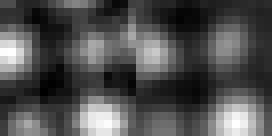}
       {BCC}
   \end{minipage}
   \begin{minipage}{0.32\linewidth}
       \centering
       \includegraphics[width=\linewidth]{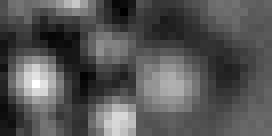}
       \includegraphics[width=\linewidth]{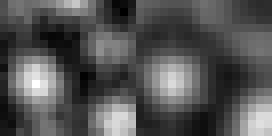}
       {FCC}
   \end{minipage}
    \caption{Comparison of original (left) and reconstructed (right) STM image patches using the More Patches configuration, showing both initial (top right) and final reconstruction (bottom right). Simple Cubic (MSE: 0.0152, SSIM: 0.9027), BCC (MSE: 0.0196, SSIM: 0.9092), and FCC (MSE: 0.0167, SSIM: 0.8768).}
   \label{fig:patch_reconstruction_comparison}
\end{figure}

The reconstructions using the More Patches configuration reveal how the model's performance varies with lattice complexity. For the Simple Cubic structure (Fig. \ref{fig:patch_reconstruction_comparison}a), the reconstruction most successfully preserves the primary structural features - the central atomic position remains well-defined, with effective noise reduction in the background regions, though at the cost of some gradient smoothing around the central feature.

The BCC structure reconstruction (Fig. \ref{fig:patch_reconstruction_comparison}b) demonstrates increased challenges with structural complexity. While atomic positions are maintained, there is more pronounced diffusion of intensity gradients in the peripheral regions. The background smoothing is more aggressive compared to the Simple Cubic case, resulting in greater loss of fine detail despite similar quantitative metrics (SSIM: 0.9092).

The FCC lattice reconstruction (Fig. \ref{fig:patch_reconstruction_comparison}c) reveals the most significant challenges. Despite reasonable quantitative performance (MSE: 0.0167), the visual quality shows notable limitations. The reconstruction exhibits the strongest background smoothing of all three types, and while the central atomic position is preserved, there is substantial loss of the surrounding structural details that characterize the FCC lattice. The gradients between atomic positions appear more homogenized, potentially obscuring important structural information.

A consistent pattern emerges across all three lattice types: as structural complexity increases from Simple Cubic to BCC to FCC, the model shows an increasing tendency to over-smooth features and lose fine structural details. This progressive degradation in reconstruction quality suggests that while our model architecture can effectively handle simpler lattice structures, it may require refinement to better preserve the more subtle structural features present in complex crystal systems. The discrepancy between quantitative metrics and visual quality, particularly evident in the FCC case, indicates that additional evaluation criteria specifically targeted at structural preservation might be valuable for future development.

\subsubsection{Performance Across Lattice Types}
Our model showed varying performance across different lattice types:
\begin{itemize}
   \item \textbf{Simple Cubic:} The model performed best on Simple Cubic structures, achieving an MSE of 0.0152 and SSIM of 0.9027 with the More Patches configuration. The reconstructions maintained good atomic position fidelity and achieved effective noise reduction while preserving essential structural features. The final reconstructed patches showed clear improvement from their initial state, particularly in background smoothing and feature definition.
   
   \item \textbf{BCC:} Performance on BCC structures showed moderate degradation compared to Simple Cubic, with best metrics of MSE 0.0196 and SSIM 0.9092 under More Patches configuration. The reconstructions exhibited increased diffusion of intensity gradients and more aggressive background smoothing, suggesting difficulty in capturing the additional complexity of BCC structures. The progression from initial to final reconstruction showed less dramatic improvement compared to Simple Cubic cases.
   
   \item \textbf{FCC:} FCC structures presented the greatest challenge, achieving MSE 0.0167 and SSIM 0.8768 with More Patches configuration. While the model maintained central atomic positions, it struggled to preserve the complex surrounding structural details characteristic of FCC lattices. The reconstructions showed the strongest tendency toward over-smoothing, with minimal improvement between initial and final states in preserving fine structural features.
\end{itemize}

A clear trend emerged across lattice types: reconstruction quality progressively degraded with increasing structural complexity. While Simple Cubic reconstructions maintained good fidelity to the original patches, BCC and FCC reconstructions showed increasing levels of detail loss and feature smoothing. This pattern suggests that our current architecture may be better suited to simpler crystal structures, with additional refinement needed for more complex lattice types.

\subsubsection{Latent Space Representation}
Analysis of the 10-dimensional latent space through PCA visualization revealed a relatively homogeneous distribution of latent vectors, without clear clustering based on lattice structures. This lack of distinct clustering suggests that the model may not be learning clearly differentiated representations for different lattice types, despite showing varying performance across them. The homogeneous distribution might indicate that the model is focusing more on local feature reconstruction rather than learning globally distinctive characteristics of each lattice type. This observation aligns with the visual reconstruction results, where the model shows stronger performance in preserving local atomic positions but struggles with larger-scale structural patterns, particularly in more complex lattice types.

\subsection{Challenges and Future Work}

Throughout this study, we encountered several challenges and identified areas for future improvement:
\begin{itemize}
    \item \textbf{Latent Space Interpretability:} The homogeneous distribution in our latent space visualization suggests that our model may not be effectively distinguishing between different lattice structures. Future work could focus on developing more interpretable latent representations, possibly through the use of disentangled variational autoencoders or by incorporating domain knowledge into the encoding process.
    
    \item \textbf{Model Complexity vs. Performance:} While our relatively shallow network (CAE-A) showed promising results, the more complex CAE-B architecture developed in MATLAB warrants further investigation. Future studies could explore the trade-offs between model complexity, computational resources, and performance gains across different lattice types.
    
    \item \textbf{Transfer Learning:} Although not implemented in the current study due to technical challenges, transfer learning presents a promising avenue for future work. This approach could potentially leverage knowledge gained from one lattice structure to improve performance on others, especially in cases where data for certain lattice types is limited.
    
    \item \textbf{Full Image Reconstruction:} Our current approach focuses on patch-level reconstruction. Extending this to full image reconstruction while maintaining structural coherence across patches is an important next step for practical STM image analysis.
    
    \item \textbf{Noise Robustness:} While our model was trained on simulated images with added noise, further work is needed to evaluate and improve its performance on real-world STM images with varying levels and types of noise.
    
    \item \textbf{Unseen Lattice Structures:} Testing the model's generalization capabilities on lattice structures not seen during training would be crucial for assessing its practical applicability in materials science research.
    
    \item \textbf{Multiple Decoders:} We plan to explore the use of multiple decoders in our future work. Specifically, we aim to develop two separate decoders: one focused on noise reduction and another on centering the image. This approach could potentially improve the model's ability to handle different aspects of image reconstruction and enhancement simultaneously.
\end{itemize}

Addressing these challenges and pursuing these future directions will be crucial in developing more robust and versatile autoencoder models for STM image analysis, potentially leading to new insights in materials science and nanotechnology. The implementation of multiple decoders, in particular, could significantly enhance the model's utility by tackling specific image processing tasks in parallel, thereby improving both the quality and interpretability of the reconstructed STM images.

\section{Conclusion}

In this study, we developed and evaluated convolutional autoencoder architectures for the analysis and reconstruction of Scanning Tunneling Microscopy (STM) images across various crystalline lattice structures. Our approach demonstrated the potential of deep learning techniques in extracting meaningful features from complex STM data, while also highlighting the challenges inherent in this task.

The CAE models showed varying degrees of success in reconstructing patches from different lattice types, with performance differences providing insights into the relative complexity of these structures. The use of simulated STM images with realistic noise characteristics allowed us to train models robust to common experimental artifacts, laying groundwork for future application to real-world STM data.

While our current results are promising, they also point to several avenues for future research. Improving the interpretability of the latent space representations, extending the model to full image reconstruction, and incorporating physical constraints are all exciting directions that could significantly enhance the utility of these models in materials science applications.

As the field of STM imaging continues to advance, the development of sophisticated analysis tools like the ones presented in this study will play an increasingly crucial role in unlocking the wealth of information contained in these atomic-scale images. By bridging the gap between raw STM data and meaningful structural insights, these autoencoder models have the potential to accelerate discoveries in nanotechnology, surface science, and beyond.

\bibliographystyle{alpha}
\bibliography{references}

\vskip.5in
Department of Mathematics

University of South Carolina

Columbia, SC 29208, USA
\end{document}